\documentclass[10pt,a4paper]{article}
\usepackage[latin1]{inputenc}
\usepackage{amsmath}
\usepackage{amsfonts}
\usepackage{amssymb}
\usepackage{graphicx}
\newtheorem{The}{Theorem}[section]
\newtheorem{Def}{Definition}[section]

\newtheorem{Ex}{Example}[section]

\begin{document}
\title{Analyzing stability of a delay differential equation involving two delays} 



\author{Sachin Bhalekar\\
Department of Mathematics,\\
 Shivaji University, Kolhapur - 416004, India\\
 Email: sbb\_maths@unishivaji.ac.in, sachin.math@yahoo.co.in}

\maketitle
\begin{abstract}
Analysis of the systems involving delay is a popular topic among applied scientists.  In the present work, we analyze the generalized equation $D^{\alpha} x(t) = g\left(x(t-\tau_1), x(t-\tau_2)\right)$ involving two delays viz. $\tau_1\geq 0$ and $\tau_2\geq 0$. We use the the stability conditions to propose the critical values of delays. Using examples, we show that the chaotic oscillations are observed in the unstable region only. We also propose a numerical scheme to solve such equations.
 
\end{abstract}


\section{Introduction}
Equations involving nonlocal operators are crucial in modeling natural systems. The memory involved in such systems cannot be modeled properly by using local operators such as ordinary integer order derivative. The fractional order derivative and the time delay (lag) are proved suitable for this task.
\par If the order of the derivative involved in the model is non-integer then it can be called as a fractional derivative. Analysis and applications of fractional derivatives can be found in \cite{SAM, POD, SHRI, FM, MAG, LFIT}.
\par If the modeling differential equation includes the past values of state variable then it is called as a delay differential equation (DDE). Basic analysis and various applications of DDE are discussed in \cite{6-SMITH, LAKSH, parkin, ECG}.
\par Stability of fractional order delay differential equations (FDDEs) is discussed by Matignon in \cite{9-MAT}. BIBO stability of linear time invariant FDDEs is discussed in \cite{PAK, BON}. A numerical algorithm is proposed by Hwang and Cheng \cite{STAB2} to test the stability of FDDEs. Stability of a class of FDDE is studied in \cite{pramana}. In  \cite{part1}, the present author proposed stability and bifurcation analysis of generalized FDDE. Chaos in various FDDEs is analyzed by Bhalekar and Daftardar-Gejji \cite{CHAOS3, CHAOS4, CHAOS5}. Application of FDDE in NMR is discussed by Bhalekar et al. \cite{BLOCH1, BLOCH2}.
\par Equations involving single delay are extensively studied in the literature. However, the analysis of the systems with multiple delay is complicated and hence the work is relatively rare. Geometry of stable regions in differential equations with two delays is presented in \cite{HALE} and references cited therein. Hopf bifurcation in these systems is studied in \cite{BELA, LIR, WUWA} and references cited therein. Applications of systems involving multiple delay can be found in various natural systems. Bi and Ruan \cite{BIR} discussed these systems in tumor and immune system interaction models. Prey-predator system with multiple delay is analyzed by Gakkhar and Singh \cite{GAKK}.

\par In this work, we extend our analysis in \cite{part1} to consider two delay differential equation with fractional order. We present stability and bifurcation analysis of these equations, propose numerical method and solve examples.

\par The paper is organized as follows: Basic definitions are listed in Section \ref{prel}. The expressions for critical curves are derived in Section \ref{mrs}. Section \ref{num} deals with the numerical method for solving general two delay models. Illustrated examples are presented in Section \ref{ille}. Conclusions are summarized in Section \ref{conc}.

\section{Preliminaries}\label{prel}
In this section, we discuss some basic definitions \cite{POD,SAM,YLG}.
\begin{Def} A real function $f(t),\, t > 0$ is said to be in space $C_\alpha,\, \alpha \in \Re$ if there exists a real number $p$ ($>\alpha$), such that $f(t)= t^pf_1(t)$ where $f_1(t)\in C[0,\infty)$.
 \end{Def}
\begin{Def}A real function $f(t),$ $t > 0$ is said to be in space $C_{\alpha}^{m}, \,  m\in I\!\!N \bigcup \left\{0\right\}$ if $f^{(m)} \in C_\alpha.$
\end{Def}
\begin{Def}Let $f \in C_\alpha$ and $\alpha \geq -1,$ then the (left-sided) Riemann--Liouville integral of order $\mu,\, \mu>0$ is given by
\begin{equation}
	I^\mu f(t)=\frac{1}{\Gamma(\mu)}\int_0^t(t-\tau)^{\mu-1}f(\tau)\,d\tau, \quad t>0.\label{2.1}
\end{equation}
\end{Def}

\begin{Def}The  Caputo fractional derivative of $f$, $f\in C_{-1}^{m}, m\in I\!\!N,$ is defined as:
\begin{eqnarray}
D^\mu f(t) && = \frac{d^m}{d t^m} f(t) ,\quad \quad \quad \mu=m \nonumber \\
&& = I^{m-\mu} \frac{d^m f(t)}{d t^m}, \quad m-1<\mu<m.\label{2.2}
\end{eqnarray}
Note that for $m-1 < \mu \leq m,\,\, m \in I\!\!N,$
\begin{eqnarray}
I^\mu D^\mu f(t)&&= f(t)-\sum^{m-1}_{k=0} \frac{d^k f}{d t^k}(0)\frac{t^k}{k!},  \label{2.3}\\
I^\mu t^\nu &&= \frac{\Gamma (\nu + 1)}{\Gamma (\mu + \nu + 1)} t^{\mu + \nu}. \label{2.4}
\end{eqnarray}
\end{Def}

\section{Main Results} \label{mrs}
We consider the following generalized delay differential equation (GDDE) 

\begin{equation}
D^{\alpha} x(t) = g\left(x(t-\tau_1), x(t-\tau_2)\right), \label{n3.1}
\end{equation}
where  $D^{\alpha}$ is a Caputo fractional derivative of order $\alpha \in (0, 1]$, $g$ is continuously differentiable function and $\tau_1 \geq 0$, $\tau_2 \geq 0$ are delays.
\par A number $x^*\in \mathbb{R}$ is called an equilibrium point  of (\ref{n3.1}) if
\begin{equation}
g\left(x^*, x^*\right)=0.\label{3.2}
\end{equation}

Linearization of the system (\ref{n3.1}) in the neighborhood of equilibrium point $x^*$ is given by

\begin{eqnarray}
D^{\alpha} \xi = a \xi_{\tau_1} + b \xi_{\tau_2}, \label{3.3}
\end{eqnarray}
 where $x_{\tau}(t)=x(t-\tau)$ and $a$, $b$ are partial derivatives of $g$ with respect to first and second variables evaluated at $\left(x^*, x^*\right)$ respectively.
The characteristic equation of the system is
\begin{equation}
\lambda^{\alpha} = a \exp(-\lambda \tau_1) + b \exp(-\lambda \tau_2).\label{3.4}
\end{equation}

\subsection{Stability of equilibrium}
We present the stability analysis by using the same strategy as in \cite{part1}. 
If all the eigenvalues $\lambda_i$ of characteristic equation (\ref{3.4}) satisfy
\begin{equation}
Re(\lambda_i)<0,\, \forall i.\label{3.4.1}
\end{equation}
then the equilibrium  $x^*$ is asymptotically stable.

Hence, the equation without delay (i.e. with $\tau_1=\tau_2=0$) is asymptotically stable if
\begin{equation}
a+b<0. \label{3.5}
\end{equation}
Now, we assmue that $\tau_1>0$, $\tau_2>0$ and  $\lambda= u+\imath v$, $u, v\in \Re$.
The stability properties of equilibrium will change at $\lambda= \imath v$. In this case, the characteristic equation takes the form
\begin{equation}
\left(\imath v\right)^\alpha = a \exp(-\imath v \tau_1) + b \exp(-\imath v \tau_2).\label{3.7}
\end{equation}
Simplifying, we get

\begin{eqnarray}
v^\alpha \cos\left(\frac{\alpha \pi}{2}\right) - a \cos(v\tau_1)  &=&   b \cos(v\tau_2) \label{3.8}\\
v^\alpha \sin\left(\frac{\alpha \pi}{2}\right) + a \sin(v\tau_1) &=&  - b \sin(v\tau_2).\label{3.9}
\end{eqnarray}

Squaring and adding
\begin{equation}
v^{2\alpha} + a^2 - 2 a v^\alpha \cos\left(\frac{\alpha \pi}{2} + v \tau_1\right)   =  b^2.
\end{equation}
This gives
\begin{equation}
\tau_1 = \frac{1}{v} \left(-\frac{\alpha \pi}{2}+\arccos\left[\frac{v^{2\alpha} + a^2-b^2}{2 a v^\alpha}\right]\right). \label{tau1}
\end{equation}
Using similar arguments and rewriting system (\ref{3.8})--(\ref{3.9}), we get
\begin{equation}
\tau_2 = \frac{1}{v} \left(-\frac{\alpha \pi}{2}+\arccos\left[\frac{v^{2\alpha} - a^2+b^2}{2 b v^\alpha}\right]\right).\label{tau2}
\end{equation}

Equations (\ref{tau1}) and (\ref{tau2}) are parametric expressions of $\tau_1$ and $\tau_2$ in parameter $v$. These can be used to plot critical curves.

\section{Numerical Method} \label{num}
In this section, we present a numerical scheme based on \cite{DSB} to solve GDDE (\ref{n3.1}) with initial function $x(t)=\phi(t),$ $t\leq 0$. This initial value problem (IVP) is equivalent to the integral equation
\begin{eqnarray}
x(t)&=&\int_{0}^{t}\frac{\left(t-\zeta\right)^{\alpha-1}}{\Gamma(\alpha)}f\left(x(\zeta-\tau_1), x(\zeta-\tau_2)\right)d\zeta\nonumber\\
&& + \phi(0), \quad t\in[0, T]. \label{eqint}
\end{eqnarray}
We consider the uniform grid $\left\{t_n=nh: n=-k, -k+1, \cdots, -1,0,1, \cdots, N\right\}$ where $k$ and $N$ are the integers chosen so that $h=T/N = \tau_1/k_1 = \tau_2/k_2$ and $k=max\left\{k_1,k_2\right\}$. Let $x_j$ be the approximation to $x(t_j)$. Note that $x_j=\phi(t_j),$ if $j\leq 0$ and $x(t_j-\tau_1)=x_{j-k_1}$, $x(t_j-\tau_2)=x_{j-k_2}$. 
\par Evaluating (\ref{eqint}) at $t=t_{n+1}$ and using product trapezoidal formula, we get
\begin{equation}
x_{n+1}=\phi(0) + \frac{h^\alpha}{\Gamma(\alpha+2)}\sum_{j=0}^{n+1}a_{j,n+1}f\left(x_{j-k_1}, x_{j-k_2}\right),
\end{equation}
where
$$
a_{j,n+1}=\left\{\begin{array}{ccc}
n^{\alpha+1}-(n-\alpha)(n+1)^{\alpha}&\mbox{ if}&j=0;\\

(n-j+2)^{\alpha+1}+(n-j)^{\alpha+1}&\\
-2(n-j+1)^{\alpha+1}&\mbox{if}&1\leq j \leq n;\\
1&\mbox{if}&j=n+1.\\
\end{array}\right.
$$

\section{Illustrative Examples}\label{ille}
In this section, we analyze stability and solve some examples using the numerical scheme presented in preceding section.

\begin{Ex}
Consider the fractional order generalization of U\c{c}ar system \cite{ucarincom}
\begin{equation}
D^\alpha x(t)=\delta x(t-\tau_1) - \epsilon \left[x(t-\tau_2)\right]^3, \label{3.1}
\end{equation}
where $\delta$ and $\epsilon$ are positive parameters.
\end{Ex}
There are three equilibrium points of the system (\ref{3.1}) viz. $0, \pm\sqrt{\delta/\epsilon}$. If $\tau_1=\tau_2=0$ then $0$ is unstable whereas $\pm\sqrt{\delta/\epsilon}$ are stable. We take $\delta=1$.
\par The characteristic equation of the system (\ref{3.1}) at $x^*= \pm\sqrt{\delta/\epsilon}$ is
\begin{equation}
\lambda^\alpha-\exp(-\lambda\tau_1)+3\exp(-\lambda\tau_2)=0.
\end{equation}
Note that the characteristic equation and hence stability does not depend on $\epsilon$. We take $\epsilon=1$. 
\par Using (\ref{tau1}) and (\ref{tau2}) we obtain the critical values of delay as

\begin{equation}
\tau_1 = \frac{1}{v} \left(-\frac{\alpha \pi}{2}+\arccos\left[\frac{v^{2\alpha} -8}{2  v^\alpha}\right]\right) \label{t1}
\end{equation}
and
\begin{equation}
\tau_2 = \frac{1}{v} \left(-\frac{\alpha \pi}{2}+\arccos\left[\frac{v^{2\alpha} + 8}{-6 v^\alpha}\right]\right). \label{t2}
\end{equation}
\par The following stability result based on \cite{LIR} for integer order case $\alpha=1$ is proposed in \cite{ucarincom}.
\begin{The}\label{Thm3} \cite{ucarincom}
The system $\dot{x}(t)=\delta x(t-\tau_1) - \epsilon \left[x(t-\tau_2)\right]^3$ is locally stable if one of the following conditions hold:\\
(i) $\tau_2\in [0, 0.366667/\delta)$\\
(ii) $\tau_2\in [0.366667/\delta, 0.43521/\delta)$ and $\tau_1\in [0, 0.0697127/\delta)$.
\end{The}
It can be verified that the stability region described in this Theorem can be readily obtained by using expressions (\ref{t1}) and (\ref{t2}). Further, as shown in Fig. 1 we can extend this region for higher values of delay also. 

Fig. 2 shows the critical curves for different values of $\alpha$.  

 \begin{figure}
 \includegraphics[scale=.7]{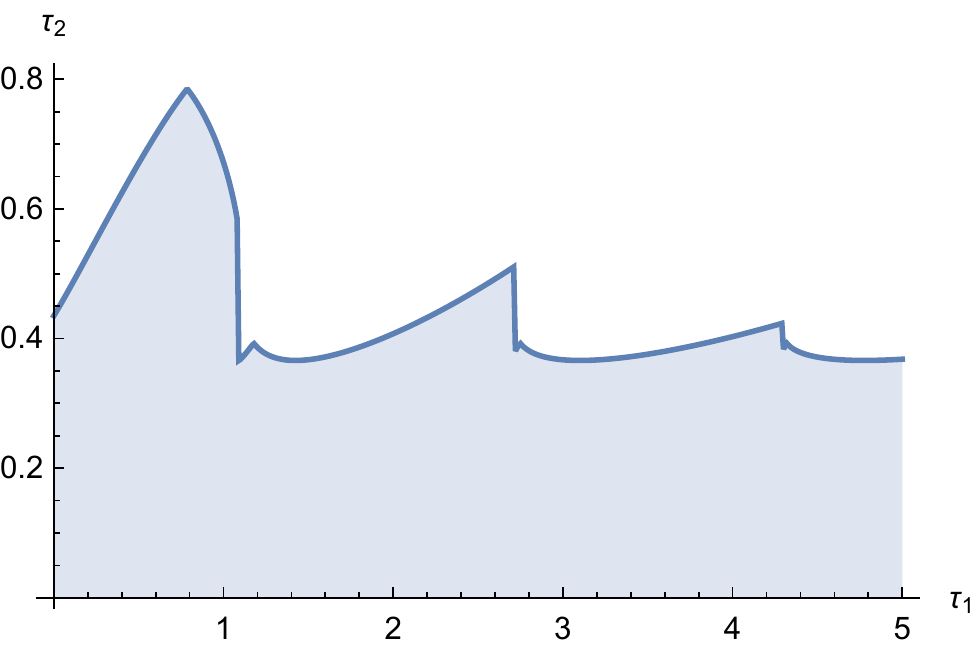}
  \caption{Extended stability region of U\c{c}ar system for $\alpha=1$ \label{Fig. 1}}
  \end{figure}
  
\begin{figure}
\includegraphics[scale=.7]{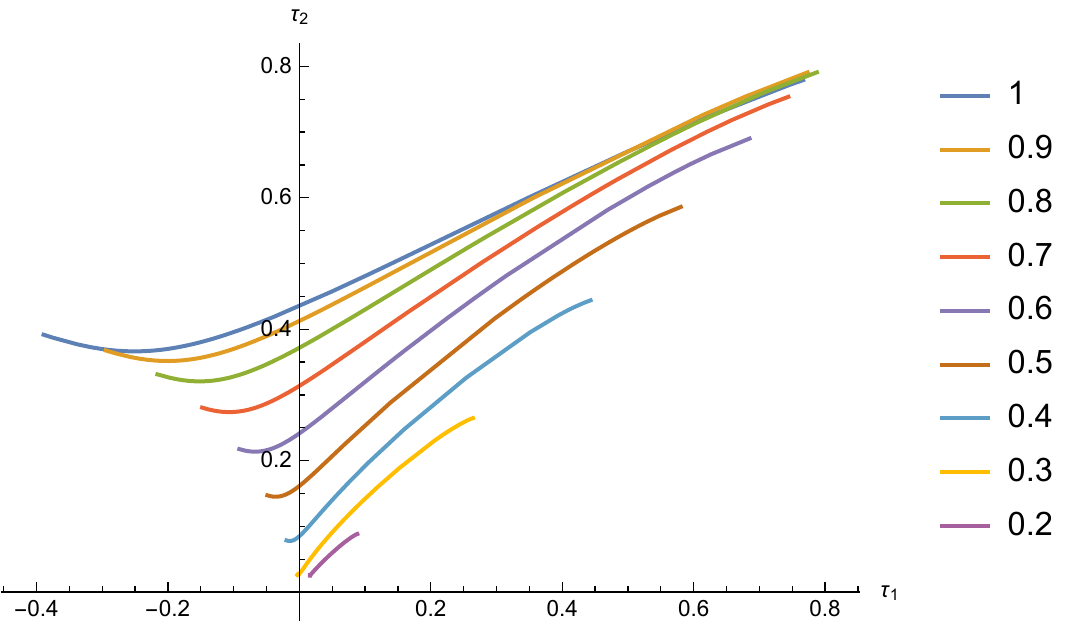}
 \caption{Critical curves of U\c{c}ar system for different values of $\alpha$ \label{Fig.  2}}
 \end{figure}
 
Now, we verify these stability results using numerical computations. In Fig. 3, we have taken $\tau_1=0.4$, $\tau_2=0.5$ and $\alpha=0.9$. These values are in stable region and the figure shows the orbit converging to an equilibrium state. In Fig. 4, we take the parameter values $\tau_1=1.6$, $\tau_2=1.4$ and $\alpha=0.9$ in unstable region. In this case, the system exhibit chaotic oscillations.

\begin{figure}
\includegraphics[scale=.7]{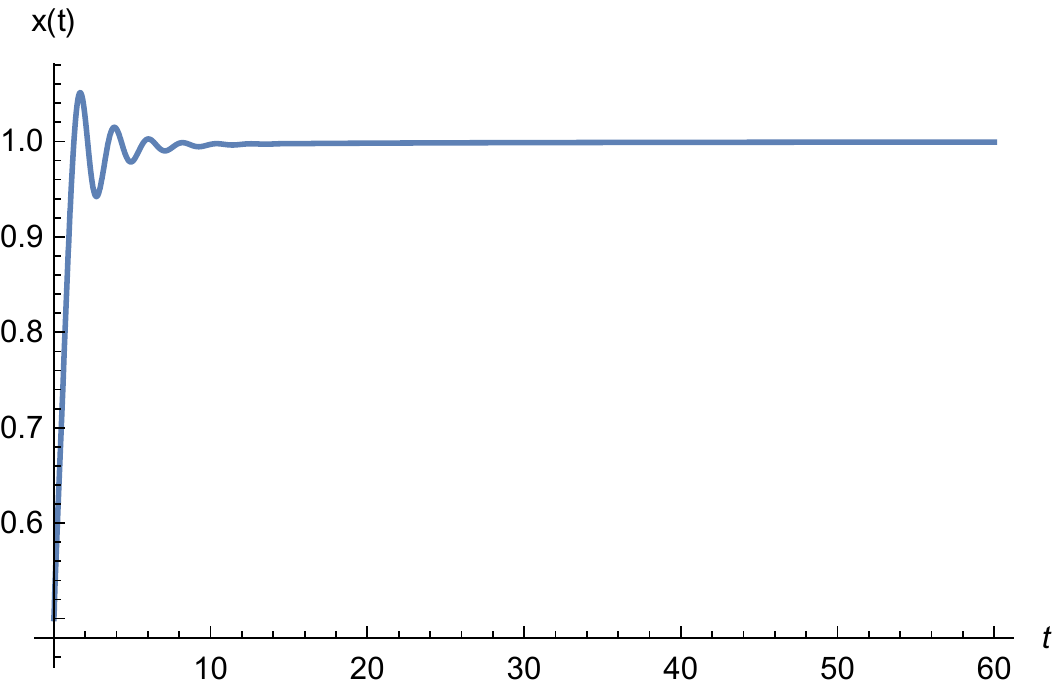}
 \caption{Stable orbit of U\c{c}ar system for $\alpha=0.9$, $\tau_1=0.4$ and $\tau_2=0.5$ \label{Fig.  3}}
 \end{figure}
 
 \begin{figure}
 \includegraphics[scale=.7]{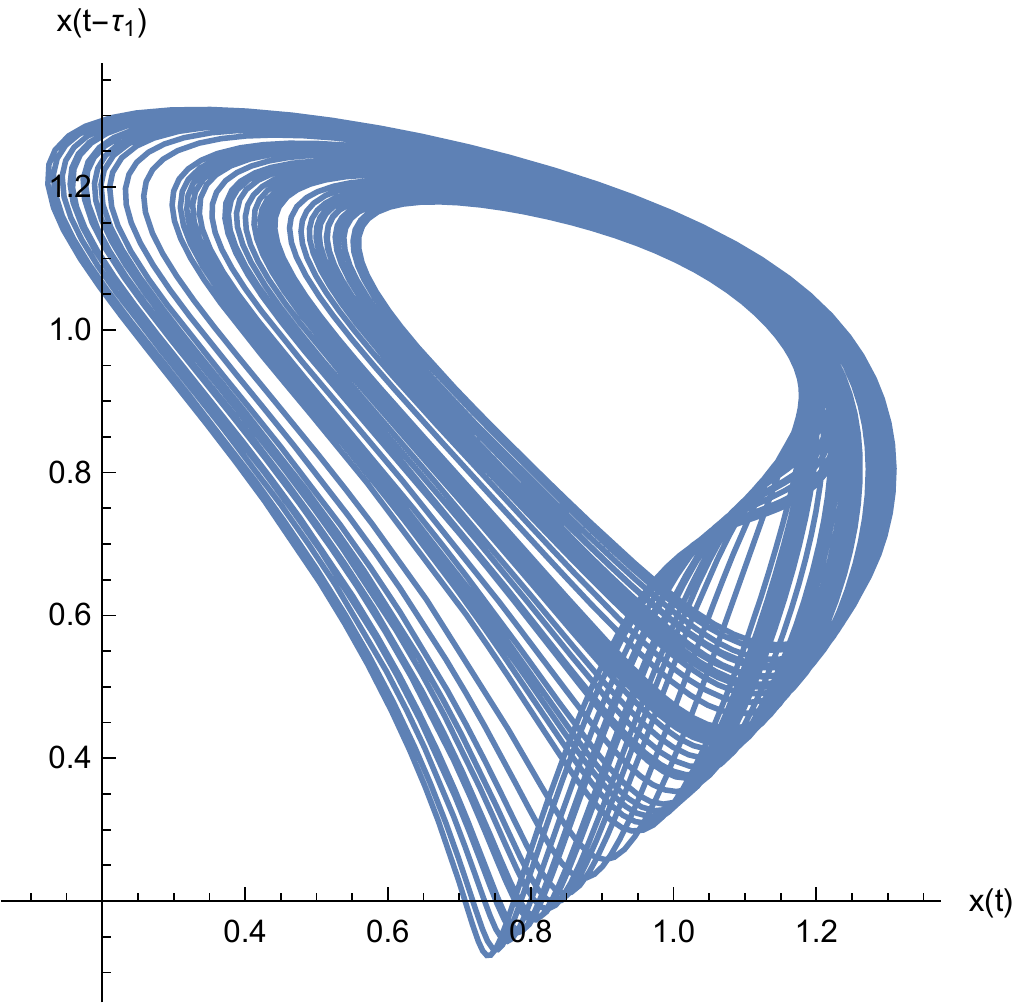}
  \caption{Chaotic attractor of U\c{c}ar system for $\alpha=0.9$, $\tau_1=1.6$ and $\tau_2=1.4$ \label{Fig.  4}}
  \end{figure}

\begin{Ex}
Consider the generalization of fractional order Ikeda equation \cite{CHAOS2} with two delays
\begin{equation}
D^\alpha x(t) = -3 x\left(t-\tau_1\right) + 24 \sin\left(x\left(t-\tau_2\right)\right),\quad 0<\alpha\leq 1. \label{ikeda}
\end{equation}
\end{Ex}
As discussed in \cite{part1}, there are seven equilibrium points. The points  $x_0=0$ and $x_{3,4}=\pm 7.49775$ are unstable whereas $x_{1,2}=\pm 2.7859$ and $x_{5,6}=\pm 7.95732$ are stable at $\tau_1=\tau_2=0$. The characteristic equation of the system (\ref{ikeda}) at the equilibrium point $x^*$ is
\begin{equation}
\lambda^\alpha=-3\exp(-\lambda \tau_1)+24 \cos(x^*) \exp(-\lambda \tau_2).
\end{equation}
If we take $x^*=x_{1,2}$ then the characteristic equation becomes
\begin{equation}
\lambda^\alpha=-3\exp(-\lambda \tau_1)-22.4977 \exp(-\lambda \tau_2).
\end{equation}
This gives $a=-3$ and $b=-22.4977$.  The critical curve in this case is given in Fig. 5 for $\alpha=0.7$.
 \begin{figure}
 \includegraphics[scale=.7]{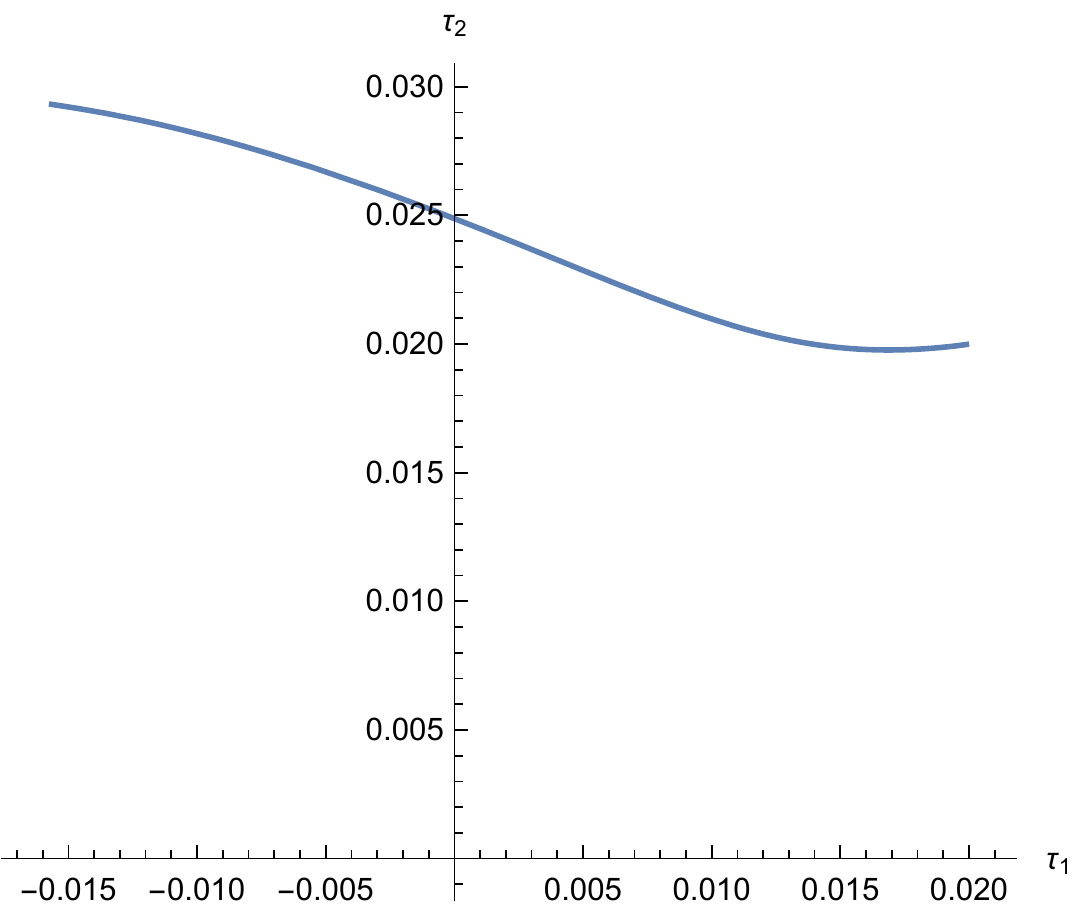}
  \caption{Critical curve of Ikeda system for $\alpha=0.7$ \label{Fig. 5}}
  \end{figure}
In Fig. 6, we have presented stable solution of system (\ref{ikeda}) with $\tau_1=0.02$ and $\tau_2=0.01$ in stable region. If we consider the values $\tau_1=0.01$ and $\tau_2=0.1$ in unstable region then we get chaotic attractor (cf. Fig. 7).

\begin{figure}
 \includegraphics[scale=.7]{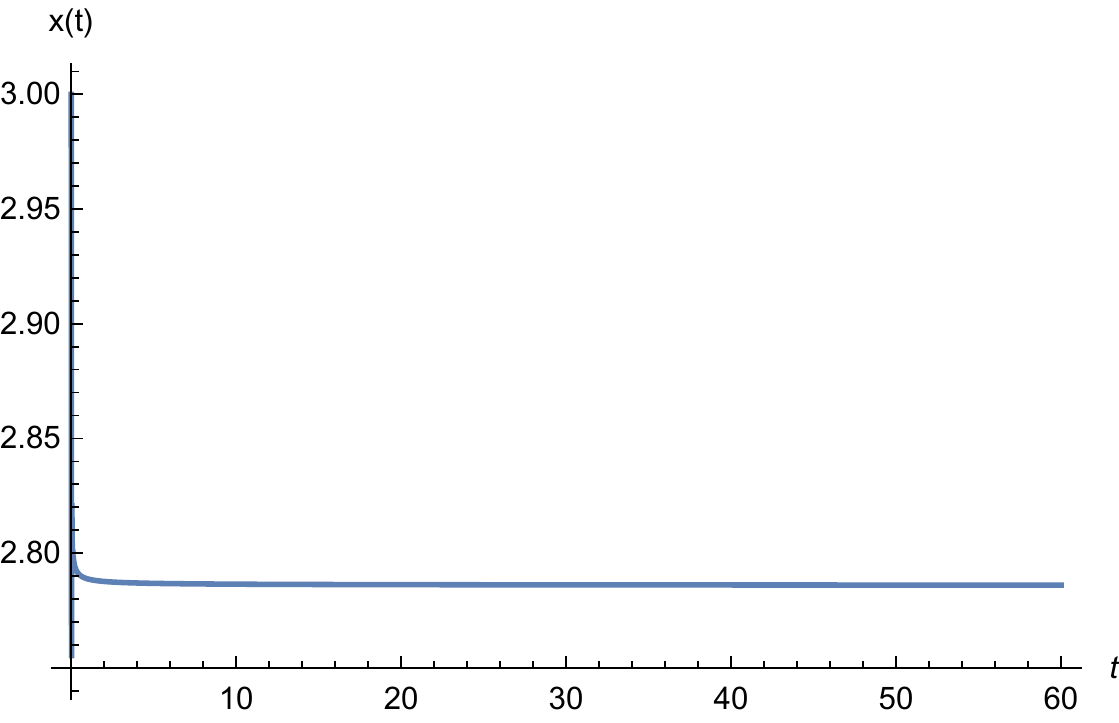}
  \caption{Stable orbit of Ikeda system for $\alpha=0.7$, $\tau_1=0.02$ and $\tau_2=0.01$  \label{Fig. 6}}
  \end{figure}
  
  \begin{figure}
   \includegraphics[scale=.7]{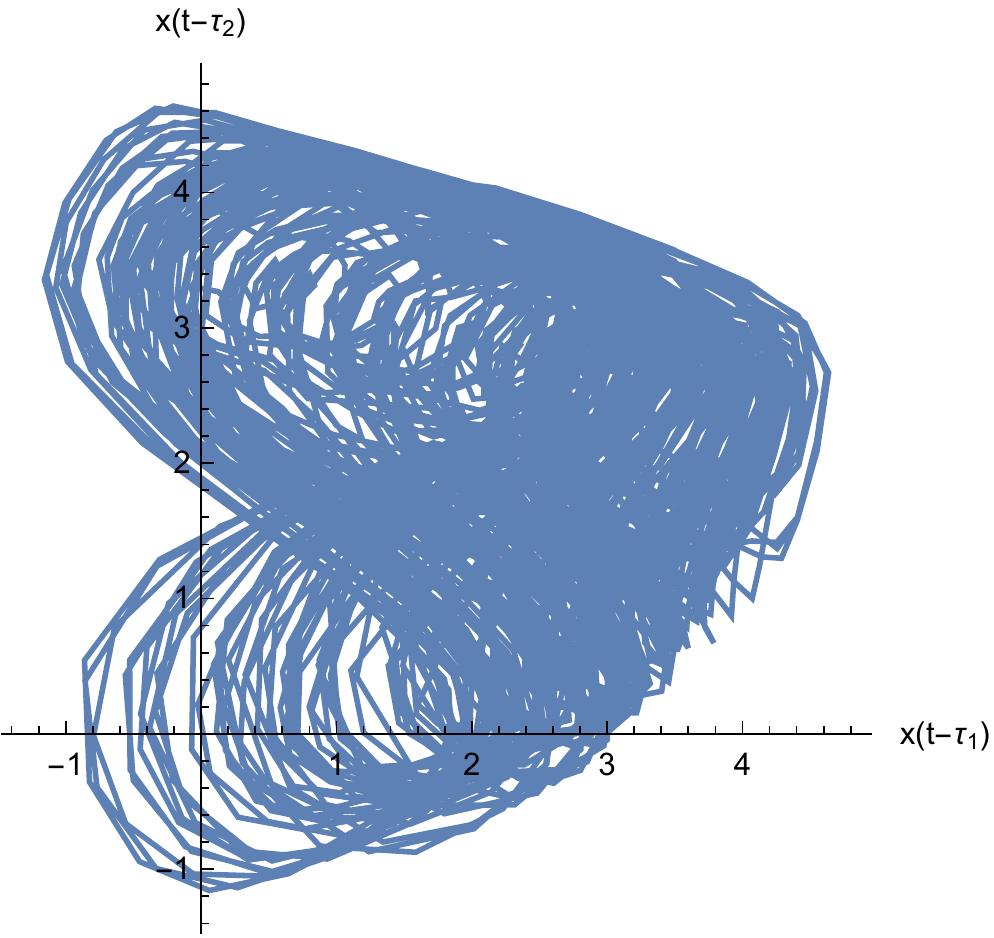}
    \caption{Chaotic attractor shown by Ikeda system for $\alpha=0.7$, $\tau_1=0.01$ and $\tau_2=0.1$ \label{Fig. 7}}
    \end{figure}

\section{Conclusion}\label{conc}
In \cite{part1}, we have presented the complete analysis of the generalized equation involving delay. However, the analysis of the systems with multiple delay is not as simple as that of single delay. 
In this article, we have presented the expressions for critical values of generalized two delay system. We obtain the parametric relation between the delays $\tau_1$ and $\tau_2$. The region bounded by horizontal axis and the critical curve is stable region if the eigenvalues are stable at $\tau_1=\tau_2=0$. The nonlinear system may exhibit chaotic oscillations in unstable region. We also have presented numerical scheme to solve these systems. Further, the results are illustrated with two examples.

\section*{Acknowledgements}
Author acknowledges the Science and Engineering Research Board (SERB), New Delhi, India for the Research Grant (Ref. MTR/2017/000068) under Mathematical Research Impact Centric Support
(MATRICS) Scheme.


\end{document}